\documentclass[mathpazo]{article}
\usepackage{amssymb,amsmath,amsthm,latexsym,graphicx}
\usepackage{fancyhdr}                

\usepackage{hyperref}
\hypersetup{linktocpage}
\hypersetup{
    colorlinks,
    citecolor=black,
    filecolor=black,
    linkcolor=black,
    urlcolor=black
}

\theoremstyle{remark}
\newtheorem{remark}{Remark}[section]

\numberwithin{equation}{section}
\numberwithin{figure}{section}

\def\df{\dfrac}

\def\tilde{\widetilde}

\def\ud{\,\mathrm{d}}
\begin{document}
\title{Numerical Resolution near $t=0$ of
Nonlinear Evolution Equations
in the Presence of Corner Singularities
in Space Dimension 1}

\author{ Qingshan Chen,
      Zhen Qin, and Roger Temam\thanks{Corresponding author: Roger
        Temam ({\tt temam@indiana.edu})}}
\date{ }
\maketitle
\thispagestyle{fancy}
\renewcommand{\headrulewidth}{0pt}
\fancyhead{ }
\fancyhead[LO]{\textsf{First published in {\em Comm.~Comp.~Phys.},
    9(3):568-586, March 2011.\\ Published by Global Science Press.
    DOI: 10.4208/cicp.110909.160310s}} 
\begin{center}
Dedicated to the memory of David Gottlieb
\end{center}

\begin{abstract}
The incompatibilities between the initial and boundary data
will cause singularities at the time--space corners, which 
in turn adversely affect the accuracy of the numerical schemes
used to compute the solutions. We study the corner singularity 
issue for nonlinear evolution equations in 1D, and propose
two remedy procedures that effectively recover much of the 
accuracy of the numerical scheme in use. Applications
of the remedy procedures to the 1D viscous Burgers equation,
and to the 1D nonlinear reaction--diffusion equation 
are presented. The remedy procedures are applicable 
to other nonlinear diffusion equations as well. 
\end{abstract}
%

%
%

%

\section{Introduction}\label{s1}
It is well known in the mathematics community that
smooth boundary and initial data do not
guarantee smooth solutions to the initial-boundary
value problems of time-dependent PDEs.
Even if the existence and uniqueness of the solution
are proved
and all the given data are as smooth as desired,
then, in order for the solution to be smooth near $t=0$,
it is necessary and sufficient that the boundary data
and the initial data satisfy an infinite set of
so-called compatibility conditions.
See
\cite{Fr64,LaSoUr68,La52, La54,
RaMa74, Sm80, Te82, Te06},
and below.

When the boundary and initial data fail to satisfy some of
the compatibility conditions, singularities may occur
 at the corner of
the time and spatial axes.
In some cases, for simple problems, this is not
a critical issue, because
the singularity is
short-lived.

However in recent years the issue of the compatibility conditions for the initial
boundary value problems of time-dependent PDEs started to
receive attention in the numerical simulation community,
because larger and more complex problems are handled thanks
to the ever--growing computing power that is available.
This produces a need to better
understand what happens during
the short initial period for certain physical processes.
Boyd and Flyer \cite{BoFl99}
analyzed the connection between incompatibility and the
rate of convergence of Chebyshev spectral series, and discussed the remedy
procedure of smoothing the initial conditions. Flyer and Swarztrauber
\cite{FlSw02}
studied the effect of the incompatibilities on the convergence rate
of spectral and finite difference methods.
Flyer and Fornberg \cite{FlFo03} proposed
a remedy procedure, based on the idea of singular corner functions, for
the heat equation, and for variable coefficient convection--diffusion
equations as well.
Flyer and Fornberg \cite{FlFo04} studied the corner basis functions for some
dispersive equations.
Bieniasz \cite{Bi05} modified the remedy procedure proposed
by Flyer and Fornberg \cite{FlFo03}, and
applied it to a diffusion-reaction system
arising from electrochemistry.

In this article we study, from a numerical point of view,
the compatibility issue for nonlinear diffusive PDEs.
We shall first present our approach in details for the
classical viscous Burgers equation in dimension one.
However our approach does not depend on any particular
property of the Burgers equation other than its diffusiveness.
Hence we believe that it applies to other nonlinear
diffusive equations as well. To demonstrate this,
we shall also apply the approach to a 1D nonlinear
reaction--diffusion equation.

For the Burgers equation
the correction procedure proposed in
\cite{FlFo03} cannot be expected to work to its full
strength, for two reasons. The first is that, for a
nonlinear PDE, the singular corner functions, which
satisfy the equation and also display corner singularities,
are hard, if not impossible, to find. The second reason
is that the superposition property does not hold for
nonlinear equations, which means
that even if in some rare cases we find the singular corner functions
for the nonlinear equation,
we cannot separate them from the solution without
introducing singular terms into the equation.
However, we
speculate here that the singular corner functions, derived
from the linearized Burgers equation, can be employed
to remove the zero order incompatibility
(see Section \ref{s2.2}). What is
less obvious is that the singular corner functions
can also be employed to remove higher order singularities
(see Section \ref{s2.3}).
However, as has been said, we cannot avoid introducing
singular terms into the equation.
To overcome the difficulty associated with
singular terms in the equation we choose the Galerkin finite
element method (FEM) as the means for constructing
the numerical scheme for the equation, because the Galerkin
FEM is based on the weak formulation of the PDE, and hence
is potentially more tolerant of the singularities in the
equation. The numerical results confirm the effectiveness
of the correction procedures that we propose.

Flyer et.~al.~\cite{FlFo03}, Bieniasz \cite{Bi05} and
the current study all seek the solutions of the target
equations (a different equation for each study) in the
form
\begin{equation}
   u = v + S,
   \label{e1.1}
\end{equation}
where $S$ is a linear combination of the singular corner
functions of the linear diffusive equations.
The major difference between our approach and
the approaches of Flyer et.~al.~in \cite{FlFo03}
and Bieniasz in \cite{Bi05} lies in the way
the correction procedure is implemented.
Both Flyer et.~al.~and Bieniasz derived the differential equation
for the new unknown $v$ and constructed the numerical
scheme for this equation directly. We instead combine
the correction procedure with the appropriate
numerical method, the Galerkin FEM in this study,
and work with the weak formulation
of the equation. In addition, Bieniasz considered only
the zeroth order incompatibility, while our study considers
the incompatibilities up to the first order.

The rest of the article is organized as follows. In Section
\ref{s2} we recall the compatibility conditions for the
viscous Burgers equation, and describe the correction
procedures for the incompatibilities between the initial
and boundary conditions. In Section \ref{s3} we
present numerical results to verify the effectiveness
of the proposed correction procedures.
In Section \ref{s3.5} we derive the correction procedures
for the 1D nonlinear reaction--diffusion equation. Numerical
results are also presented.
We conclude with Section \ref{s4}.

All the equations considered in this article
and in the references quoted above relate to
evolution equations in space dimension one.
In an article in preparation \cite{ChQinTe2}
we will consider higher spatial dimensions which
necessitate totally different methods.
\section{The numerical scheme and correction procedures }\label{s2}
\subsection{ The Viscous Burgers equation}\label{s2.1}
We consider the initial and boundary value problem
for the Burgers equation:
\begin{equation}\left\{\label{e2.1}
   \begin{aligned}
   &u_t + uu_x - \nu u_{xx} = 0,\qquad 0 < x < 1,\, t > 0, \\
   &u(0,t) = g_1(t),\quad u(1,t) = g_2(t),\qquad t > 0,\\
   &u(x,0) = h(x),\qquad 0 < x < 1,
   \end{aligned}\right.
\end{equation}
where $\nu$ is a small positive parameter representing the viscosity,
and
$g_1$, $g_2$ and $h$ are given real functions, assumed to be
as smooth as desired. The exact problem (\ref{e2.1}) is known 
to have a unique solution for all times; this equation is 
indeed similar, but simpler, than the 2 dimensional 
incompressible Navier-Stokes Equations for which the existence 
and uniqueness of the solution is known (see e.g. \cite{Te01}). 
Looking for a weak solution, if $h$ is given in $L^{2}(0,1)$ 
and $g_{1},g_{2}$ in $H^{\frac{1}{2}}(0,T)$, then $u$ exists 
and is unique in 
$\mathcal{C}([0,T];L^{2}(0,1))\cap L^{2}(0,T;H^{1}(0,1))$. 
If, as we assume here, $h,g_{1},g_{2}$ are smooth, then $u$ is smooth up to $\mathcal{C}^{\infty}$ regularity in $[0,1]\times (0,T]$. Now the fact that the interval $(0,T]$ is open at $0$ is related to the compatibility problem we are addressing here; even if $h,g_{1},g_{2}$ are given $\mathcal{C}^{\infty}$ (in respectively $[0,1]$ and $[0,T]$). For the (unique) solution to be smooth (even $\mathcal{C}^{0},\mathcal{C}^{1},\mathcal{C}^{2}$) near $t=0$, that is in $[0,1]\times [0,T]$, $h,g_{1},g_{2}$ must satisfy certain compatibility conditions as described in the references quoted \cite{RaMa74, Sm80,Te06}. We now make explicit the first and second compatibility conditions for (\ref{e2.1}) which guarantee respectively that $u$ is $\mathcal{C}^{0},\mathcal{C}^{1}$, in $[0,1]\times[0,T]$.\\

The compatibility conditions
require that at the corners of the time and spatial axes,
the derivatives of the solutions computed through the
boundary conditions be equal to those computed through the
Cauchy-Kowalevski rules, that is, for the viscous Burgers equation
above,
\begin{align*}
&0^\mathrm{th}\textrm{ order:} & &g_1(0) = h(0),& &g_2(0) = h(1) \\
&1^\mathrm{st}\textrm{ order:}& &g_{1t}(0) = - h(0)h'(0)+\nu h''(0) ,&
&g_{2t}(0) = - h(1)h'(1)+\nu h''(1) \\
& & &\cdots\cdots&   &\cdots\cdots
\end{align*}

We note here that the compatibilities
at the left and right corners can be treated separately, and the method
presented below apply to the incompatibilities at
both corners.
For this reason, and for the sake of simplicity,
we study, in this article, the case in which
 the compatibility
conditions at $x=1$ are met, at least to
a certain order, but those at the
left $x=0$ are not. Assuming so,   we  let
\begin{align*}
&\alpha_0\equiv g_1(0) - h(0)\neq 0,\\
&\alpha_1\equiv g_{1t}(0) + h(0)h'(0)-\nu h''(0) \neq 0,\\
&\cdots\cdots
\end{align*}
Nonzero $\alpha_0$, $\alpha_1$ represent incompatibilities
at the left time--space corner.

\subsection{ Correction procedure 1 for
the zeroth order incompatibility}\label{s2.2}
We aim to separate the singular part of the solution
from the nonsingular part by using certain singular
corner functions, but this approach will inevitably
introduce singular terms into the equation,
because the singular corner functions for
the nonlinear equation $\eqref{e2.1}_1$
is hard, if not impossible, to find,
and also because the superposition
property does not hold for nonlinear
equations. We, however, speculate  that the correction
procedure can be employed to remove the zero
order incompatibility, which is the most significant
one.
To overcome the difficulty associated with
the singular terms in the equation we choose
Galerkin finite element method (FEM) as the means
for constructing the numerical scheme for \eqref{e2.1},
because Galerkin FEM is based on the weak formulation
of the PDE, and therefore is potentially more
tolerant of the singularities in the equation.

The weak formulation of $\eqref{e2.1}_1$ can be formally
derived as follows.
Let $v$ be a continuous function that vanishes at
$x=0,\,1$. Multiplying $\eqref{e2.1}_1$ by $v$ and integrating
the equation by parts over $(0,\,1)$, we obtain
\begin{equation}\label{e2.1a}
(u_t,\,v) + \df{1}{2}(u^2,\,v_x) + \nu(u_x,\,v_x) = 0,
\end{equation}
where $(f,\,g)=\int_0^1fg\ud x$ for any $L^2$ integrable
functions.

We introduce the corner function:
\begin{equation}\label{e2.2}
S_0 = \dfrac{1}{\sqrt{\pi\nu t}}\int^\infty_x e^{-\frac{s^2}{4\nu t}}\ud s =
\mathrm{erfc}(\df{x}{2\sqrt{\nu t}}),
\end{equation}
and we notice that $S_0$ satisfies the heat equation and
displays a singularity at $(0,\,0)$:
\begin{equation}\label{e2.3}
\left\{\begin{aligned}
&S_{0t} - \nu S_{0xx} = 0,\\
&S_0(0,t) = 1,\\
&S_0(x,0) = 0.
\end{aligned}\right.
\end{equation}

We let $N$ be the number of segments in the interval
$[0,\,1]$, and $h=1/N$,
and let $V_h$ be the finite element space, with
basis $\{\varphi_0,\,\varphi_1,\,\cdots,\,\varphi_N\}$.
We look for a solution of \eqref{e2.1} in the form
\begin{equation}
u(x,t) = \alpha_0S_{0}(x,t) + v_h(x,t),
\label{e2.4}
\end{equation}
where $v_h(\cdot,t)\in V_h \textrm{ for a.e. } t > 0 $.
Imposing the boundary conditions $\eqref{e2.1}_2$, and noticing
$\eqref{e2.3}_{2}$, we have
\begin{equation}
\left\{\begin{aligned}
& v_h(0,t) = g_1(t) - \alpha_0,\\
& v_h(1,t) = g_2(t) - \alpha_0S0(1,t).
\end{aligned}\right.
\label{e2.3b}
\end{equation}
Imposing the initial condition $\eqref{e2.1}_3$, and noticing
$\eqref{e2.3}_3$, we have
\begin{equation}
v_h(x,0) = h(x),\qquad 0<x<1.
\label{e2.3c}
\end{equation}
We observe that the zeroth order incompatibility at the left time--space
corner has been removed for
$v_h$; indeed
\begin{equation}\label{e2.3a}
v_h(0,0) = g1(0) - \alpha_0 = h(0).
\end{equation}
The compatibility conditions at the right corner are not affected
because $S_0$ is smooth there.

We then plug \eqref{e2.4} into \eqref{e2.1a} and take
$v=\tilde v_h\in V_h$, with $\tilde v_h(0) = \tilde v_h(1) = 0$,
and we obtain
\begin{equation}
(\alpha_0S_{0t} + v_{ht},\,\tilde v_h) -
\df{1}{2}( (\alpha_0S_0 + v_h)^2,\, \tilde v_{hx}) +
\nu (\alpha_0S_{0x} + v_{hx},\, \tilde v_{hx} ) = 0.
\label{e2.5}
\end{equation}
Noticing that $S_0$ satisfies the heat equation, we rewrite \eqref{e2.5}
as
\begin{equation}
( v_{ht},\,\tilde v_h) -
\df{1}{2}(  v_h^2,\, \tilde v_{hx}) -
(\alpha_0S_0v_h,\,\tilde v_{hx}) +
\nu ( v_{hx},\, \tilde v_{hx} ) =
\df{1}{2}\alpha_0^2(S_0^2,\,\tilde v_{hx}).
\label{e2.6}
\end{equation}
Coupling \eqref{e2.6} with the boundary conditions \eqref{e2.3b}
and the initial condition \eqref{e2.3c} we can find $v_h$.

After we find $v_h$ we recover the original solution
of \eqref{e2.1} by \eqref{e2.4}.

Numerical experiments will be presented
in Section \ref{s3}.
For the tests that we have done, the errors
during the initial period of time is reduced
by a magnitude of more than one order.

\subsection{Correction procedure 2 for the zeroth and first order
incompatibilities}\label{s2.3}
To further improve the accuracy we
consider the next (first order) incompatibility and
introduce the second corner function:
\begin{equation}\label{e2.8a}
S_1 = \int^t_0 S_0(x,\tau)\ud\tau.
\end{equation}
We notice that the function $S_1$ satisfies
\begin{equation}\label{e2.4a}
\left\{\begin{aligned}
&S_{1t} - \nu S_{1xx} = 0,\\
&S_1(0,t) = t,\\
&S_1(x,0) = 0.
\end{aligned}\right.
\end{equation}

We look for a solution of \eqref{e2.1} in the form
\begin{equation}
u(x,t) = \alpha_0S_0(x,t) + \alpha_1S_1(x,t) + v_h(x,t),
\label{e2.9}
\end{equation}
where $v_h(\cdot,t)\in V_h$ for almost every $t>0$.
Imposing the boundary conditions $\eqref{e2.1}_2$, and noticing
$\eqref{e2.3}_2$ and $\eqref{e2.4a}_{2}$, we have
\begin{equation}
\left\{\begin{aligned}
& v_h(0,t) = g_1(t) - \alpha_0 - \alpha_1t,\\
& v_h(1,t) = g_2(t) - \alpha_0S_0(1,t) - \alpha_1S_1(1,t).
\end{aligned}\right.
\label{e2.7}
\end{equation}
Imposing the initial condition $\eqref{e2.1}_3$, and noticing
$\eqref{e2.3}_3$ and $\eqref{e2.4a}_3$, we have
\begin{equation}
v_h(x,0) = h(x),\qquad 0<x<1.
\label{e2.7a}
\end{equation}
As for the zeroth order correction procedure, this
correction procedure also removes the zeroth order
incompatibility for $v_h$ at the left corner;
indeed
\begin{equation}
v_h(0,0) = g_1(0) - \alpha_0 = h(0).
\label{e2.8}
\end{equation}
The singular corner function $S_1$ has no effect on the zero
order incompatibility, but helps to remove the first order
incompatibility. To see this, we plug \eqref{e2.9} into
the original equation $\eqref{e2.1}_1$, and noticing
that $S_0$ and $S_1$ satisfy the heat equations
\eqref{e2.3} and \eqref{e2.4} respectively, we have
\begin{equation}
v_{ht} + \left( (\alpha_0S_0+\alpha_1S_1)v_h +
\df{1}{2}(\alpha_0S_0+\alpha_1S_1)^2\right)_x + v_hv_{hx}
-\nu v_{hxx} = 0.
\label{e2.9a}
\end{equation}
We first calculate $v_{ht}(0,0)$ using the boundary condition
$\eqref{e2.7}_1$,
\begin{equation}
v_{ht}(0,0) = g_{1t}(0) - \alpha_1.
\label{e2.10}
\end{equation}
Then applying Cauchy--Kowalevsky rule to equation \eqref{e2.9a}, we
have
\begin{multline}
   v_{ht}(0,0) =\\ \left.\left[ -
   \left((\alpha_0S_0+\alpha_1S_1)v_h
   +\df{1}{2}(\alpha_0S_0+\alpha_1S_1)^2\right)_x
   - v_hv_{hx}
   + \nu v_{hxx} \right]\right|_{(0,0)}.
   \label{e2.11}
\end{multline}
Noticing that $S_0(x,0) = S_1(x,0) =0$
(see $\eqref{e2.3}_3$ and $\eqref{e2.4a}_3$), we have
\begin{equation}
   v_{ht}(0,0) = -h(0)h_x(0) + \nu h_{xx}(0).
   \label{e2.12}
\end{equation}
The right hand sides of \eqref{e2.10} and \eqref{e2.12} are
equal by the definition of $\alpha_1$.

This correction procedure does remove the zeroth and first
order incompatibilities, but it introduces singular terms
into the equation. To overcome this difficulty we again
construct the numerical scheme by Galerkin FEM, for the
reason already mentioned before.
As for the zeroth order correction procedure,
 we plug \eqref{e2.9} into \eqref{e2.1a} and take
$v=\tilde v_h\in V_h$, with $\tilde v_h(0) = \tilde v_h(1) = 0$,
and we obtain
\begin{multline}
( v_{ht},\,\tilde v_h) -
\df{1}{2}(  v_h^2,\, \tilde v_{hx}) -
\left( (\alpha_0S_0+\alpha_1S_1)v_h,\,\tilde v_{hx} \right) + \\
\nu ( v_{hx},\, \tilde v_{hx} )=
\df{1}{2}\left( (\alpha_0S_0+\alpha_1S_1)^2,\,\tilde v_{hx}\right).
\label{e2.13}
\end{multline}
We supplement \eqref{e2.13} with the boundary conditions
\eqref{e2.7} and initial condition \eqref{e2.7a},
and solve the resulting system for $v_h$. Then
we recover the solution $u$ of \eqref{e2.1} by \eqref{e2.9}.

Our results showed
that, when this procedure is applied, the accuracy
of the result is further improved.
See Section \ref{s3.2}.

\section{Numerical implementation of the correction
procedures and the results}\label{s3}
\subsection{The numerical schemes}\label{s3.1}
To unify the presentation of the numerical schemes
for the different correction procedures
we introduce
\begin{equation}\label{e3.3}
S = \left\{\begin{aligned}
     &0& &\textrm{ if no correction procedure is applied},\\
     & \alpha_0S_0&
           &\textrm{ if Correction procedure  1 is applied},\\
     &\alpha_0S_0+\alpha_1S_1&
           &\textrm{ if Correction procedure 2  is applied}.
\end{aligned}\right.
\end{equation}
The boundary conditions for $v_h$, \eqref{e2.3b} or
\eqref{e2.7} depending on the correction procedure,
can be written in one common form:
\begin{equation}\left\{\label{e3.4}
   \begin{aligned}
   &v_h(0,t) = g_1(t)- S(0,t),\\
   &v_h(1,t) = g_2(t) - S(1,t).
   \end{aligned}\right.
\end{equation}
The equation for $v_h$, \eqref{e2.6} or \eqref{e2.13},
can also be written in one common form,
\begin{equation}
( v_{ht},\,\tilde v_h) -
\df{1}{2}(  v_h^2,\, \tilde v_{hx}) -
(Sv_h,\,\tilde v_{hx}) +
\nu ( v_{hx},\, \tilde v_{hx} ) =
\df{1}{2}(S^2,\,\tilde v_{hx}),
\label{e3.5}
\end{equation}
for every $\tilde v_h\in V_h$ with
$\tilde v_h(0) = \tilde v_h(1) = 0$.

The incompatibilities between
the initial and boundary conditions
have a more severe effect
on higher order schemes than on
lower order schemes (see \cite{BoFl99,FlSw02}).
For the basis of the finite element space $V_h$ (introduced
in Section \ref{s2.2}) we choose piecewise linear functions,
which usually provide  second order approximations
to smooth functions.
The effectiveness of the correction
procedures is already evident
with the resulting numerical scheme.
We leave the endeavor for higher order schemes to future
work.

Let $N$ be the number of segments in the interval $[0,1]$,
$\Delta x=\df{1}{N}$, and $x_j = j\Delta x$ for $0\leq j\leq N$.
Let $\phi_j$ be the piecewise linear hat functions:
\begin{figure}[ht]
\scalebox{0.6}{\includegraphics{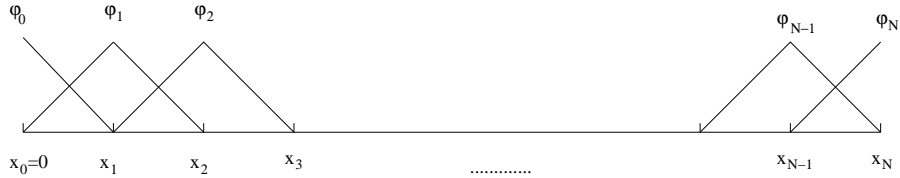}}
\caption{Piecewise linear finite elements}
\label{H}
\end{figure}
Then we write $v_h$ as a linear combination of these
basis functions:
\begin{equation}\label{e3.2}
v_h(x,t) =  \sum_{n=0}^N v_n(t)\phi_n(x),
 \end{equation}
and $v_n$, for $0\leq n\leq N$, are the unknowns.

Imposing the boundary conditions \eqref{e3.4} on
$v_h$ we obtain
\begin{equation}\left\{\label{e3.6}
   \begin{aligned}
   &v_0 = g_1(t)- S(0,t),\\
   &v_N = g_2(t) - S(1,t).
   \end{aligned}\right.
\end{equation}

Next we plug \eqref{e3.2}
into \eqref{e3.5}, and take $\tilde v_h = \varphi_m$,
with $1\leq m \leq N-1$,
\begin{multline}\label{e3.7}
\sum_{n=0}^Nv_{nt}(\phi_n,\,\phi_m) -
  \df{1}{2}\sum_{n=0}^Nv_{n}^2(\phi_n,\,\phi_{mx})
  - \sum_{n=0}^Nv_n(S\phi_n,\,\phi_{mx}) \\
  + \nu\sum_{n=0}^Nv_{n}(\phi_{nx},\,\phi_{mx}) =
  \df{1}{2}(S^2,\,\phi_{mx}),
  \textrm{ for } 1\leq m\leq N-1.
\end{multline}
This set of equations, supplemented with the
boundary conditions \eqref{e3.6} can be
easily integrated by an ODE solver.

\subsection{The results}\label{s3.2}
For demonstration purpose, we take as a test case
$\nu = 0.2$, $g_1 = 0$, $g_2 = 0$
and $h = -\sin(\df{5\pi x}{4} + \df{3\pi}{4})$.
It is easy to check that, for this test case,
both the zeroth and first order compatibility
conditions at the right corner are met, but
those at the left corner are not.

To study the accuracy of the numerical scheme we need
a means to measure the errors in the solution.
Given arbitrary initial and boundary
conditions in \eqref{e2.1}, generally no analytic
solution is known for the Burgers equation, and hence
there is no way to compute the real
errors. As an alternative, we compute
the {\it comparative} errors,
which are the differences between
two numerical solutions for the problem,
one with the stated mesh sizes,
and the other with finer mesh sizes.
In what follows the term error
is to be understood in this sense.

\begin{figure}[ht]
\scalebox{0.7}{\includegraphics{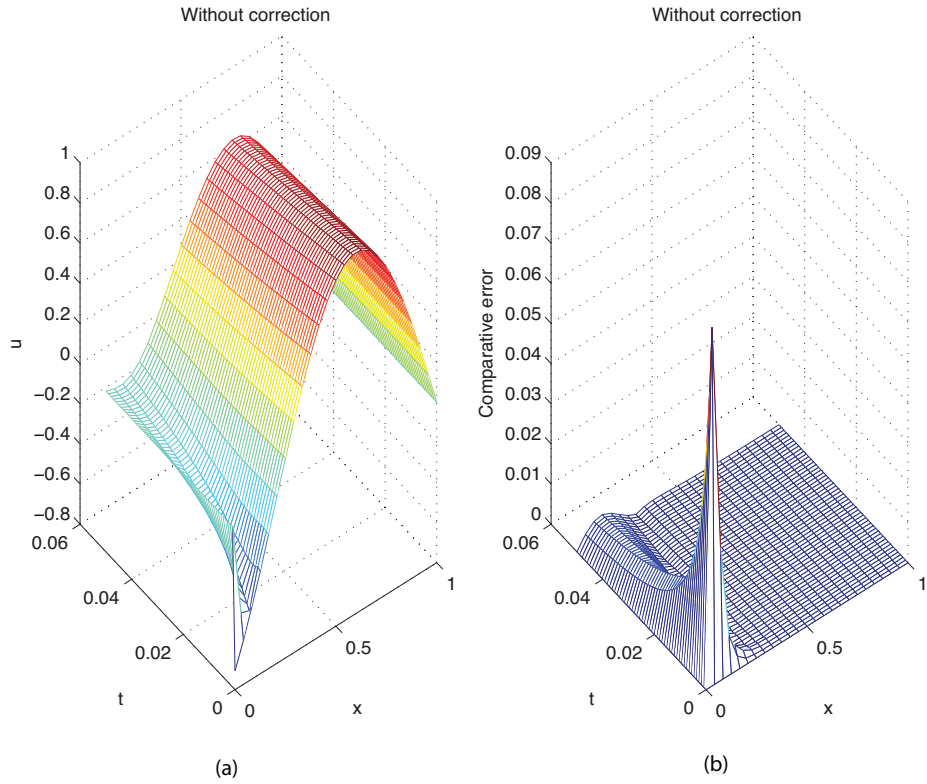}}
\caption{The solution, and pointwise errors
without any correction procedure.}
\label{f1}
\end{figure}

\begin{figure}[ht]
\scalebox{0.7}{\includegraphics{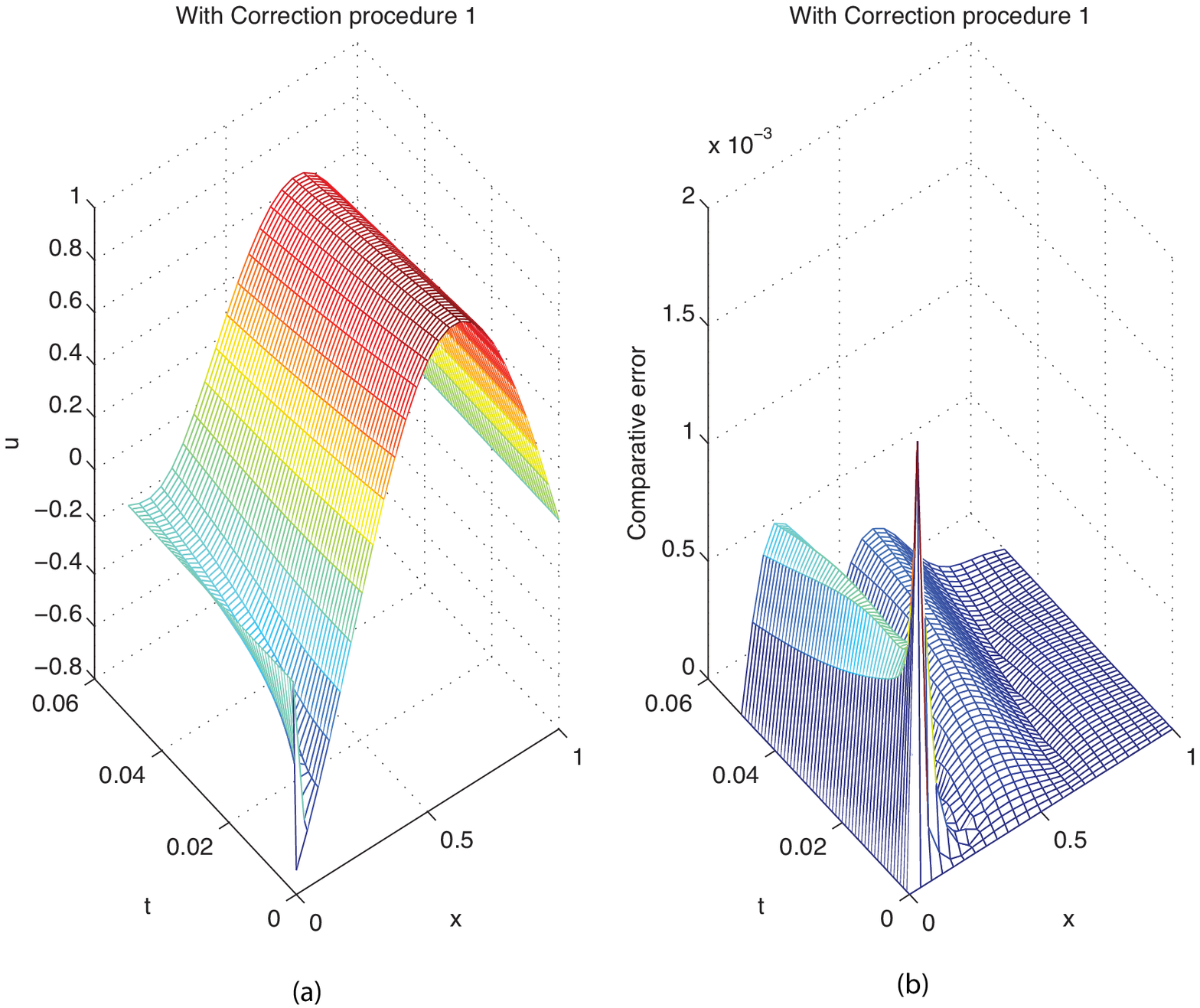}}
\caption{The solution, and pointwise errors
with {\it Correction procedure 1} applied.}
\label{f1.5}
\end{figure}

\begin{figure}[ht]
\scalebox{0.7}{\includegraphics{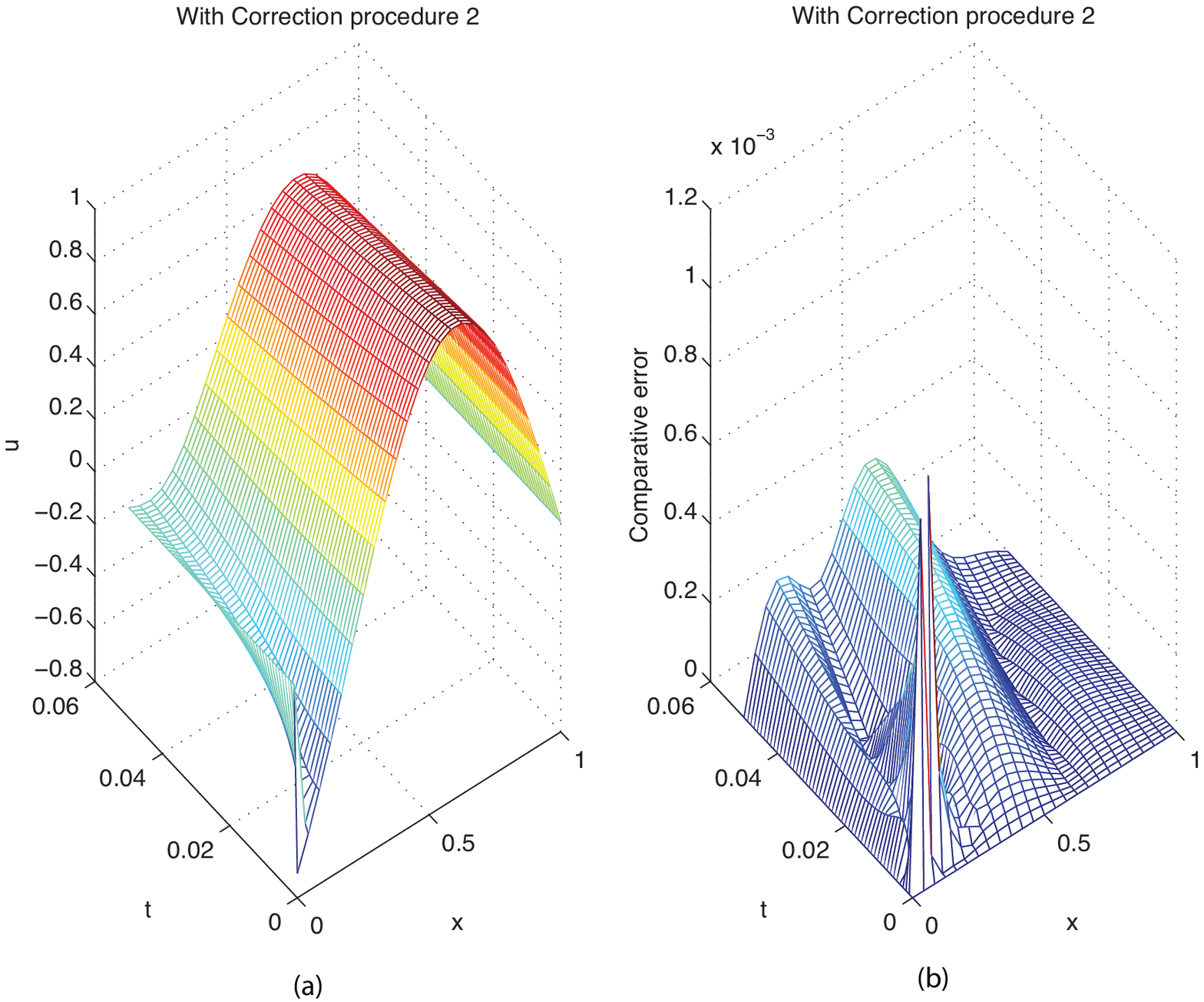}}
\caption{The solution, and pointwise errors
with {\it Correction procedure 2} applied.}
\label{f1.8}
\end{figure}

We first compute the solution of \eqref{e2.1} without any
correction procedure, i.e.~with $S=0$ in \eqref{e3.6} and
\eqref{e3.7}.
The solution is plotted in Fig.~\ref{f1} (a), and it
displays sharp gradient around the left corner
of the time--space axes due to the incompatibility
between the initial and boundary conditions there.
In order to have an overview
of the structure of the errors in the solution we plot
the pointwise errors in Fig.~\ref{f1} (b). The pike appear
near the left corner of the time--space axes, as expected,
and it dissipates away quickly.
For comparison we plot, in Fig.~\ref{f1.5}, the solution and the pointwise
errors computed with {\it Correction Procedure 1},
and, in Fig.~\ref{f1.8}, the solution and the pointwise errors
computed with {\it Correction procedure 2}.
With {\it Correction procedure 1}
the magnitude of the errors at the left corner (see Fig.~\ref{f1.5} (b))
has been reduced by roughly two orders;
with {\it Correction procedure 2}
the magnitude of the errors at the left corner (see Fig.~\ref{f1.8} (b))
are further reduced.

\begin{figure}[ht]
\scalebox{0.8}{\includegraphics{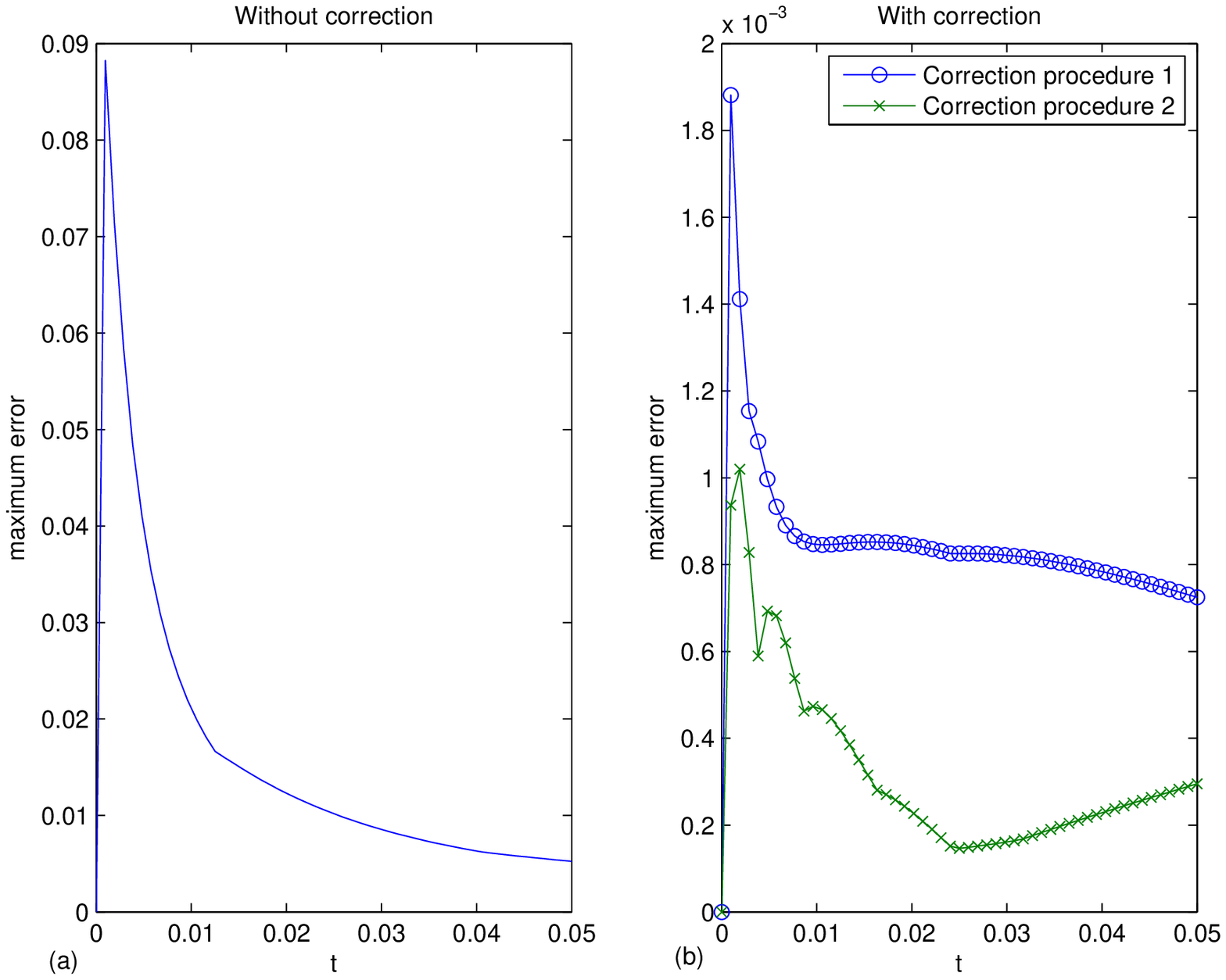}}
\caption{ Evolution of the maximum errors till $t=0.05$.
(a) Without correction procedure. (b) With correction
procedure 1 and 2.}
\label{f3}
\end{figure}

The evolution of the maximum errors at each time step
is plotted in Fig.~\ref{f3} (a).
The evolution of the maximum errors, with
{\it Correction Procedure 1}
 enabled, is displayed in Fig.~\ref{f3} (b).
The magnitude of the maximum errors has been
reduced by roughly two orders
(compared to Fig.~\ref{f3} (a)).

When \emph{Correction Procedure 2} is enabled,
we see further improvement in the accuracy, though
it is less dramatic. For comparison we plot the result
in the same figure as that for the result with
\emph{Correction Procedure 1}, and we see that the
magnitude of the maximum errors is roughly halved.

\begin{remark}\label{r1}
In practice a correction procedure should be enabled
during a short period of time at the beginning,
and be disabled afterwards, when
the solution has become smooth enough, to avoid
the computational burden associated with the singular
terms in \eqref{e3.5}. For computations that produced
Fig.~\ref{f3}, however, we run the simulation
for a short initial period of time only, and enable
the correction procedure for the whole process to
avoid unnecessary complications in programming.
\end{remark}
\begin{remark}\label{r2}
The zone of rapid variation of the solution and of 
the error should not be confused with a boundary 
layer effect. Indeed, firstly $\nu=0.2$ is too large 
to produce a sharp boundary layer, as the boundary 
layer size is of order $\sqrt{\nu}=0.447$ which is 
essentially of order 1. Furthermore when $\nu$ is 
small enough to produce a sharp boundary layer, it 
usually appears {\it for all times}, at either or 
both ends of the interval, $x=0,1$, whereas in our 
example the zone of rapid variation is concentrated 
near $x=0$, for a short time. Regarding the boundary 
layers for Burgers equation, see e.g. \cite{AN,ZWK}, 
and also \cite{CTMK}. We intend to address in a 
future work the issue of the incompatible data as 
a boundary layer at time $t=0$; see as already a first 
result in this direction in \cite{ChQinTe2}.
\end{remark}

\subsection{Comparison of convergence rates}\label{s3.3}
In this subsection we study the decay of the maximum errors,
with and without the correction procedures applied.
The results also demonstrate the effectiveness of
the correction procedures.

\begin{figure}[ht]
\scalebox{0.7}{\includegraphics{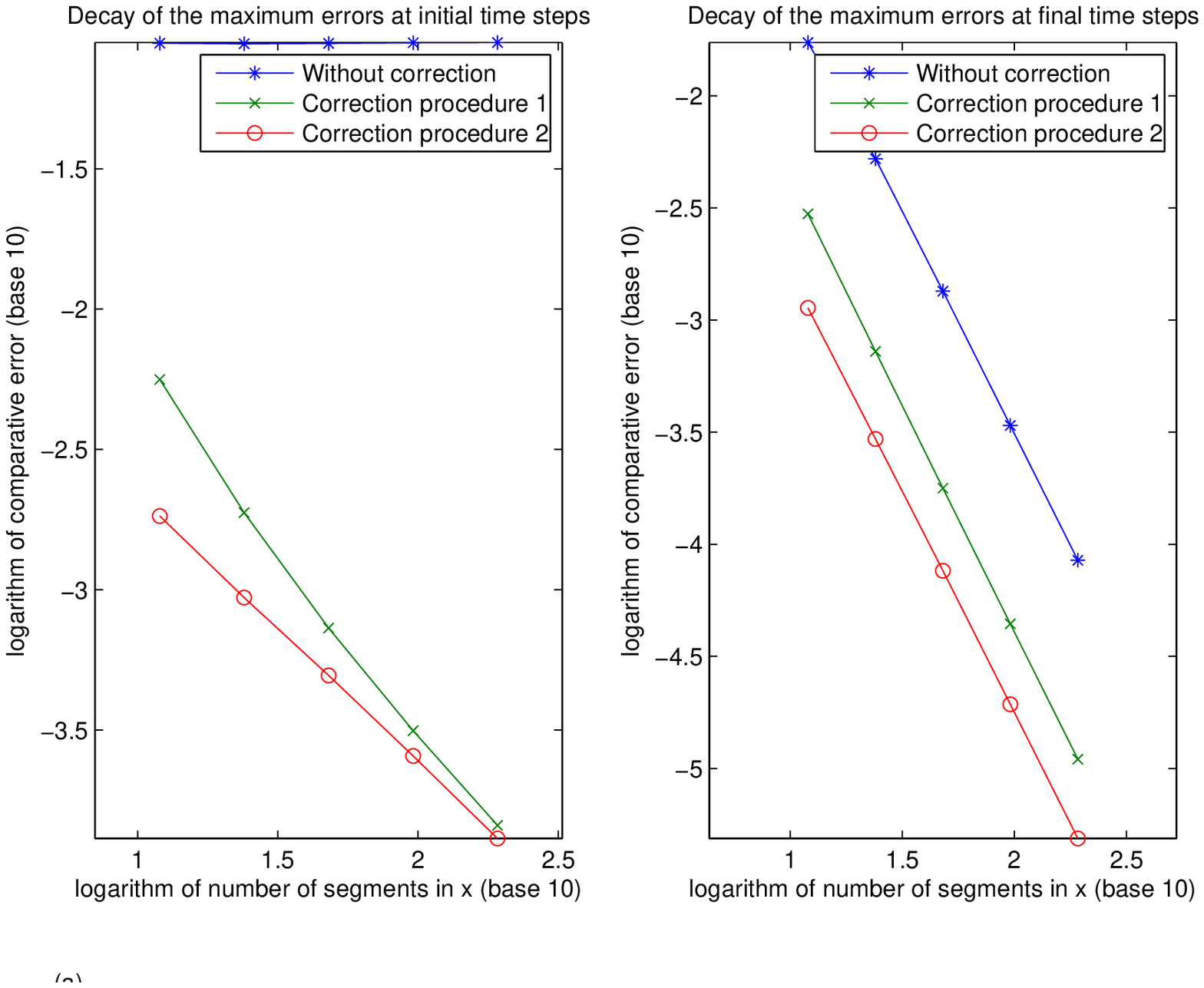}}
\caption{ Decay of the maximum errors,
(a) at the initial time step , (b) at final time step
($t=0.05$).}
\label{f4}
\end{figure}
Fig.~\ref{f4} (b) shows that, with and without the
correction procedure applied,  the maximum errors
at a fixed time $t=0.05$ decay at approximately
the second order (the slope of the line),
which is the order of accuracy of the finite element
scheme. However, the maximum errors of the simulation
with Correction Procedure 1 applied is at a magnitude
about one order smaller than without any correction procedure.
The maximum errors of the simulation with Correction Procedure
2 applied is smaller by another magnitude of about $0.3$.

The most interesting and informative comparison
can be made between the decay rates of the maximum
errors at the initial time step ($t=\Delta t$, $\Delta t$
varying with different configurations).
In Fig.~\ref{f4} (a),
the curve for the simulation without any correction
procedure stick to the upper frame of the figure;
the maximum errors wont come down whatsoever.
With {\it Correction procedure 1}, the maximum errors
decay at roughly the first order with respect
to $\Delta x$, and with {\it Correction procedure 2}, the
results are slightly better in terms of magnitude
of the maximum errors.

\section{A nonlinear reaction-diffusion equation}\label{s3.5}
In this section we extend the previous correction procedures to the nonlinear reaction-diffusion equation:\\
\begin{equation}\left\{\label{e3.1}
\begin{aligned}
&u_{t}-\nu u_{xx}+p(u)=0, 0<x<1, t>0,\\
&u(0,t)=g_{1}(t),u(1,t)=g_{2}(t), t>0,\\
&u(x,0)=h(x),
\end{aligned}\right.
\end{equation}\\
where $\nu>0$ is a parameter representing viscosity;
$g_{1}(t),g_{2}(t)$ and $h(x)$ are given functions;
and $p(u)$ is a polynomial in $u$ of odd order and
with a positive leading coefficient.
The effectiveness of these correction procedures will
be demonstrated by numerical results.

We first derive the compatibility conditions between
the initial and boundary conditions of \eqref{e3.1},
as explained in Section \ref{s2.1}:
\begin{equation}\label{e3.1a}
0^{th} order:\hspace{0.4cm} \qquad g_{1}(0)=h(0),\qquad \qquad \qquad \hspace{0.21cm}g_{2}(0)=h(1),
\end{equation}
\begin{equation}\label{e3.1b}
1^{st} order:g_{1t}(0) = \nu h''(0)-p(h(0)) ,\quad g_{2t}(0) = \nu h''(1)-p(h(1)).
\end{equation}
For simplicity, in what follows, we only consider the case where incompatibilities are present only at the left corner. Thus we let\\
\begin{align}
   &\alpha_{0}\equiv g_{1}(0)-h(0)\neq 0,\label{e3.1c}\\
   &\alpha_{1}\equiv g_{1t}(0)-\nu h''(0)+p(h(0))\neq 0.\label{e3.1d}\\
&\cdots\cdots   \qquad\cdots\cdots
\end{align}
Of course the method we present for deriving the correction
procedures would also apply to the incompatibilities at the right corner.

\subsection{The correction procedure}\label{s3.51}
As we did for the nonlinear convection--diffusion equation in Section \ref{s2} and \ref{s3}, we shall combine the correction procedures with the weak formulation of $(\ref{e3.1})_{1}$, for the reason stated in Section \ref{s2.2}. We multiply $(\ref{e3.1})_{1}$ by $v\in \mathfrak{D}(0,1)$ and integrate by parts to obtain the weak formulation of the equation:\\
\begin{equation}\label{e4.1}
(u_{t},v)+\nu (u_{x},v_{x}) + (p(u),v)=0.
\end{equation}\\

\vspace{0.1cm}
We intend to present the correction procedures (corrections of the incompatibilities to various order) in a unified way. To this end we shall employ the following notation:\\
\begin{equation}\label{e4.2}
S = \left\{\begin{aligned}
     &0& &\textrm{ if no correction procedure is applied},\\
     & \alpha_0S_0&
           &\textrm{ if Correction procedure  1 is applied},\\
     &\alpha_0S_0+\alpha_1S_1&
           &\textrm{ if Correction procedure 2  is applied}.
\end{aligned}\right.
\end{equation}\\
Here the singular corner functions $S_{0}$ and $S_{1}$ are defined as in (\ref{e2.2}) and (\ref{e2.8a}) respectively. We let N be the number of segments in the interval $[0,1], h=1/N$ and let $V_{h}$ be the finite element space spanned by $(\varphi_{0},\varphi_{1},\cdots ,\varphi_{N})$ (see Fig.1).\\

\vspace{0.1cm}

We look for the solutions of (\ref{e4.1}) in the form \\
\begin{equation}\label{e4.3}
u\simeq S+v_{h}(x,t),
\end{equation}
where  $v_{h}(\cdot,t)\in V_{h}$ for a.e.~$t$. Plugging (\ref{e4.3}) into (\ref{e4.1}), and taking $v=\tilde{v_{h}}\in V_{h}$ with $\tilde{v}_{h}(0)=\tilde{v}_{h}(1)=0$, we obtain\\
\begin{equation}\label{e4.4}
(v_{ht},\tilde{v}_{h})+\nu (v_{hx},\tilde{v}_{hx}) + (p(v_{h}+S),\tilde{v}_{h})=0.
\end{equation}
Imposing the boundary conditions and initial conditions in (\ref{e3.1}) we have\\
\begin{equation}\label{e4.5}
\left\{\begin{aligned}
&v_{h}(0,t)=g_{1}(t)-S(0,t),\\
&v_{h}(1,t)=g_{2}(t)-S(1,t),
\end{aligned}\right.
\end{equation}
and\\
\begin{equation}\label{e4.6}
v_{h}(x,0)=h(x).
\end{equation}
We will solve (\ref{e4.4}), (\ref{e4.5}) and (\ref{e4.6}) for $v_{h}$, and then we recover $u$ by (\ref{e4.3}).\\

\vspace{0.1cm}

Concerning the compatibility conditions for $v_{h}$ we make the following observations. For the zeroth order correction procedure, $S=\alpha_{0}S_{0}$, and by (\ref{e4.5}) and (\ref{e2.3}), we have \\
\begin{equation}\label{e4.7}
\left\{\begin{aligned}
&v_{h}(0,t)=g_{1}(t)-\alpha_{0},\\
&v_{h}(1,t)=g_{2}(t)-\alpha_{0}S_{0}(1,t).
\end{aligned}\right.
\end{equation}
The initial and boundary conditions for $v_{h}$ satisfy the zeroth order compatibility condition. Indeed,\\
\begin{equation}\label{e4.8}
g_{1}(0)-\alpha_{0}=h(0)
\end{equation}
The compatibility conditions at the right corner are not affected.\\

\vspace{0.1cm}

For the first order correction procedure $S=\alpha_{0} S_{0}+\alpha_{1}S_{1}$, and by (\ref{e4.5}), (\ref{e2.3}) and (\ref{e2.4a}),\\
\begin{equation}\label{e4.9}
\left\{\begin{aligned}
&v_{h}(0,t)=g_{1}(t)-\alpha_{0}-\alpha_{1}t,\\
&v_{h}(1,t)=g_{2}(t)-\alpha_{0}S_{0}(1,t)-\alpha_{1}S_{1}(1,t).
\end{aligned}\right.
\end{equation}
It is easy to see that the initial and boundary condition for $v_{h}$ satisfy the zeroth order compatibility conditions. They also satisfy the first order compatibility condition. To see this, we insert (\ref{e4.3}) into $(\ref{e3.1})_{1}$ and obtain\\
\begin{equation}\label{e4.10}
v_{ht}-\nu v_{hxx}+p(\alpha_{0}S_{0}+\alpha_{1}S_{1}+v_{h})=0.
\end{equation}
\begin{figure}[ht]
\scalebox{0.8}{\includegraphics{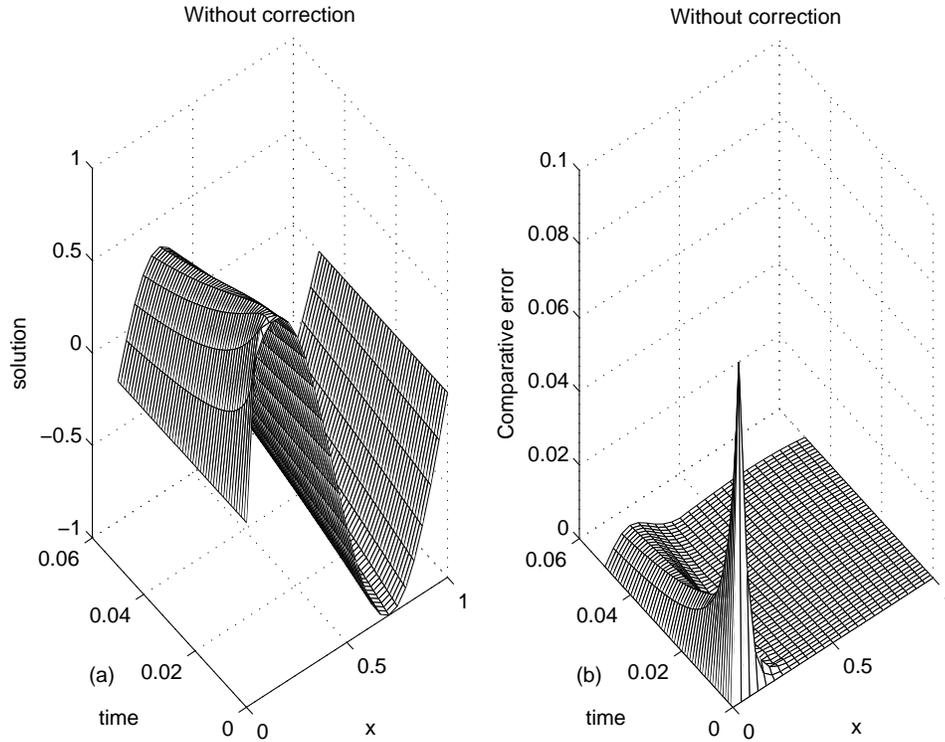}}
\caption{The solution, and pointwise errors
without any correction procedure for $p(u)=u^{3}$.}
\label{f5}
\end{figure}
We first calculate $v_{ht}(0,0)$ from(\ref{e4.9}):\\
\begin{equation}\label{e4.11}
v_{ht}(0,0)=g_{1t}(0)-\alpha_{1}
\end{equation}
Then we calculate the same quantity from (\ref{e4.10}), using the initial condition (\ref{e4.6}) instead:\\

\begin{equation}\label{e4.12}
v_{ht}(0,0)=\nu h_{xx}(0)-p(h(0)).
\end{equation}
The right-hand sides of (\ref{e4.11}) and (\ref{e4.12}) are equal according to (\ref{e3.1d})\\
\begin{figure}[ht]
\scalebox{0.8}{\includegraphics{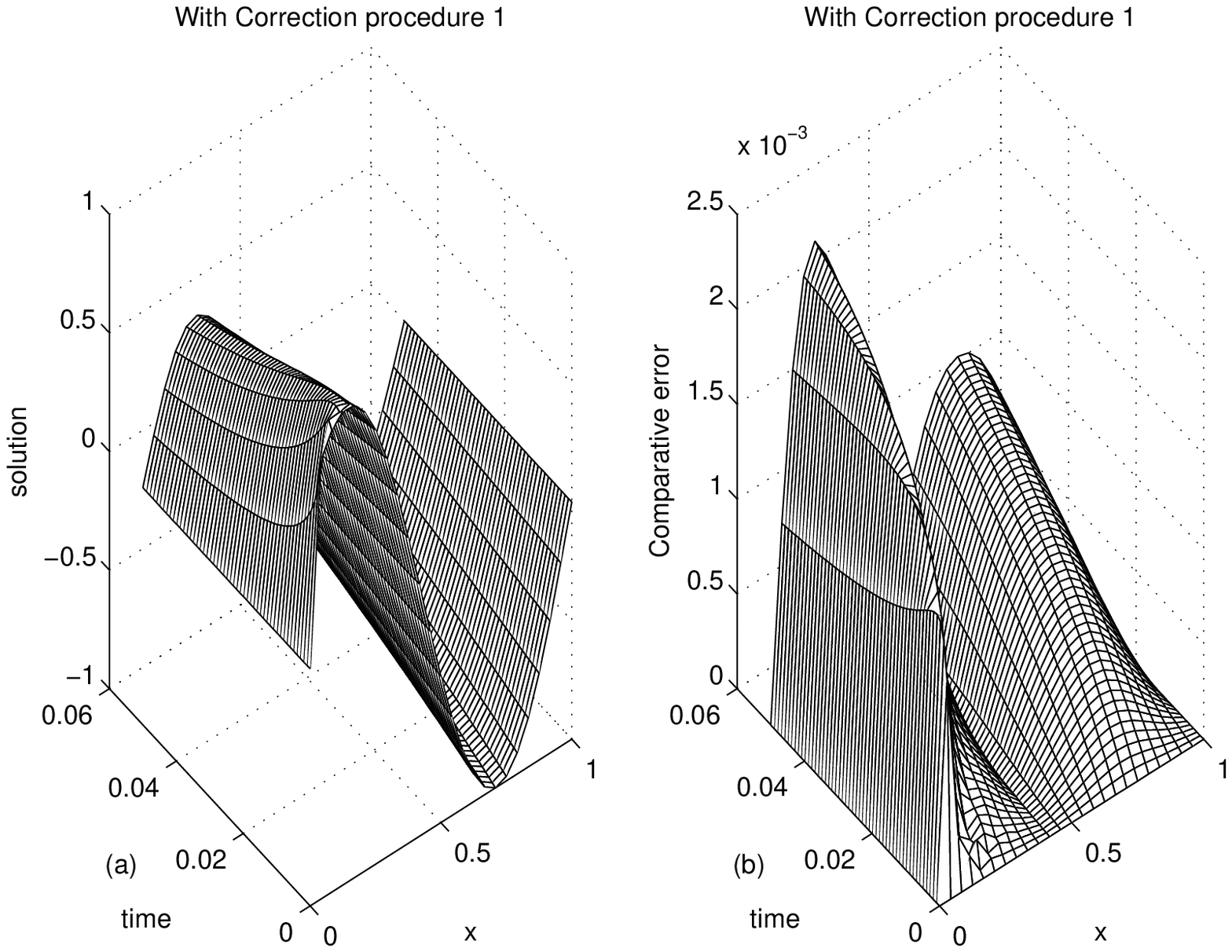}}
\caption{The solution, and the pointwise errors with Correction procedure 1 applied}
\label{f6}
\end{figure}
\begin{figure}[ht]
\scalebox{0.8}{\includegraphics{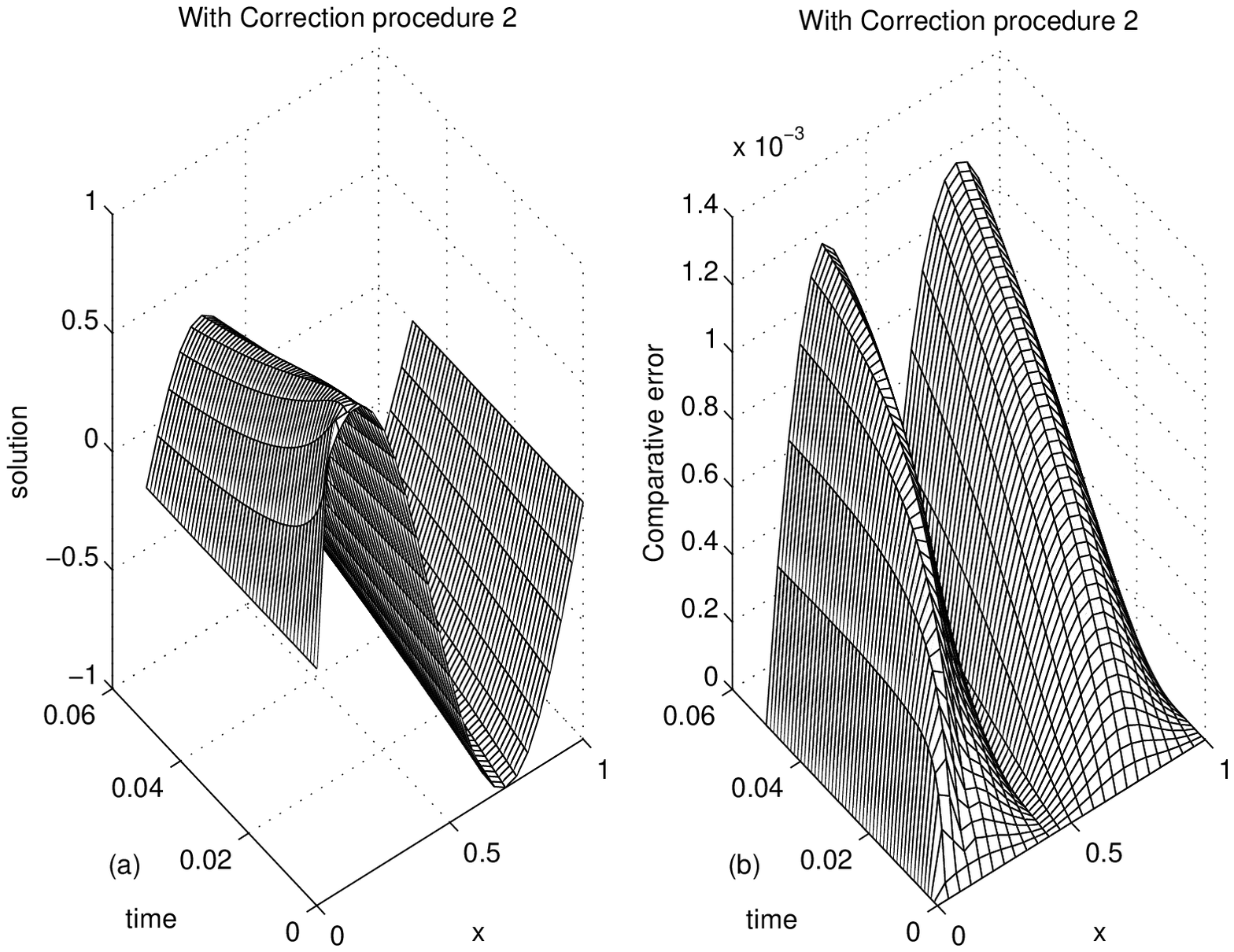}}
\caption{The solution, and the pointwise errors with Correction procedure 2 applied}
\label{f9}
\end{figure}

\begin{figure}[ht]
\scalebox{0.71}{\includegraphics{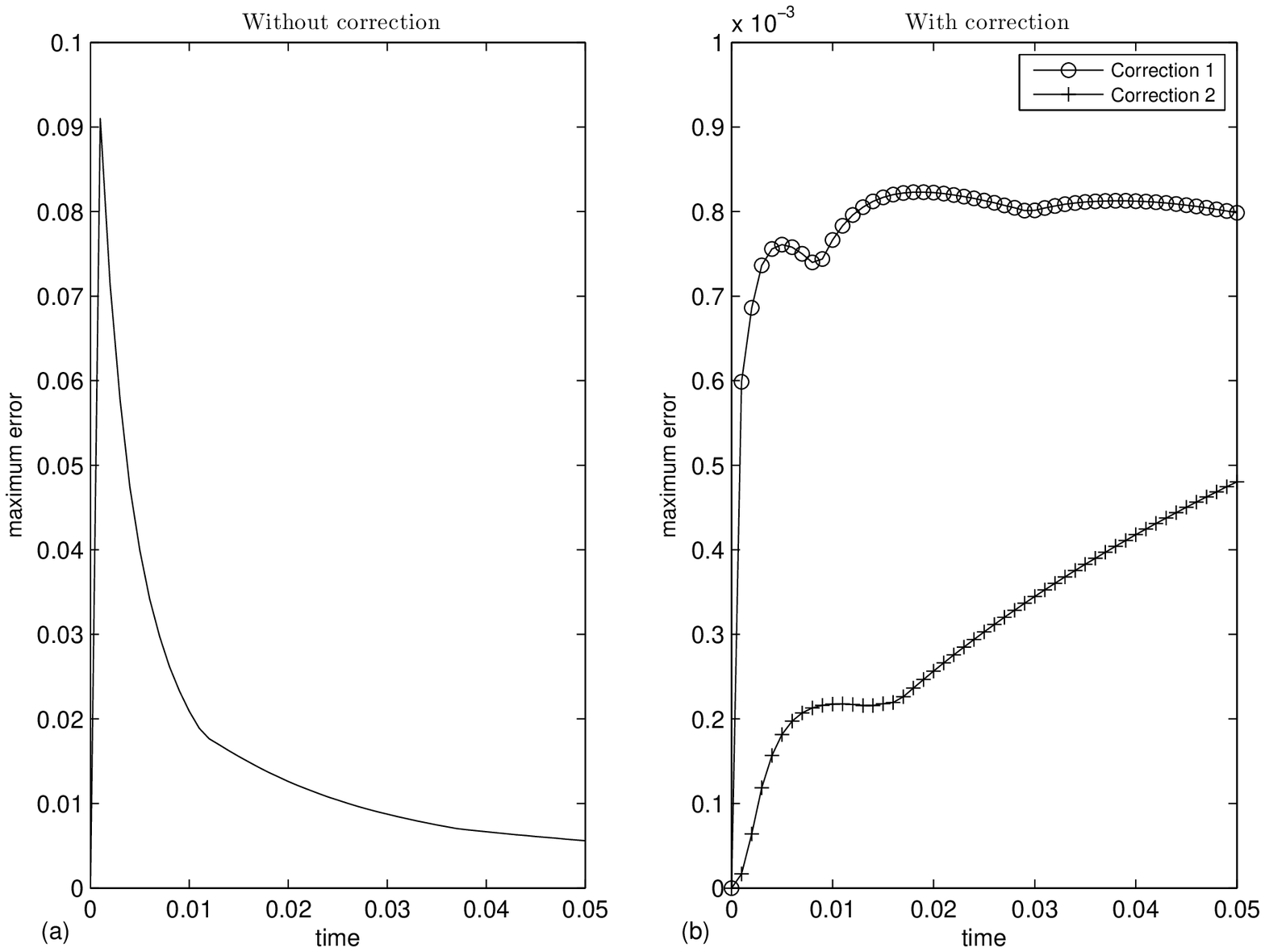}}
\caption{The evolution of the maximum errors till $t=0.05$: (a) without correction; (b) with Correction procedures 1 and 2.}
\label{f7}
\end{figure}

\subsection{The numerical results and interpretation}\label{s4.2}
In the above we presented the correction procedures for reaction-diffusion equation with polynomial reaction terms. For the numerical testing of these correction procedures we consider the specific case where\\
\begin{equation}\label{e4.13}
p(u)=u^{3}.
\end{equation}
\begin{figure}[ht]
\scalebox{0.75}{\includegraphics{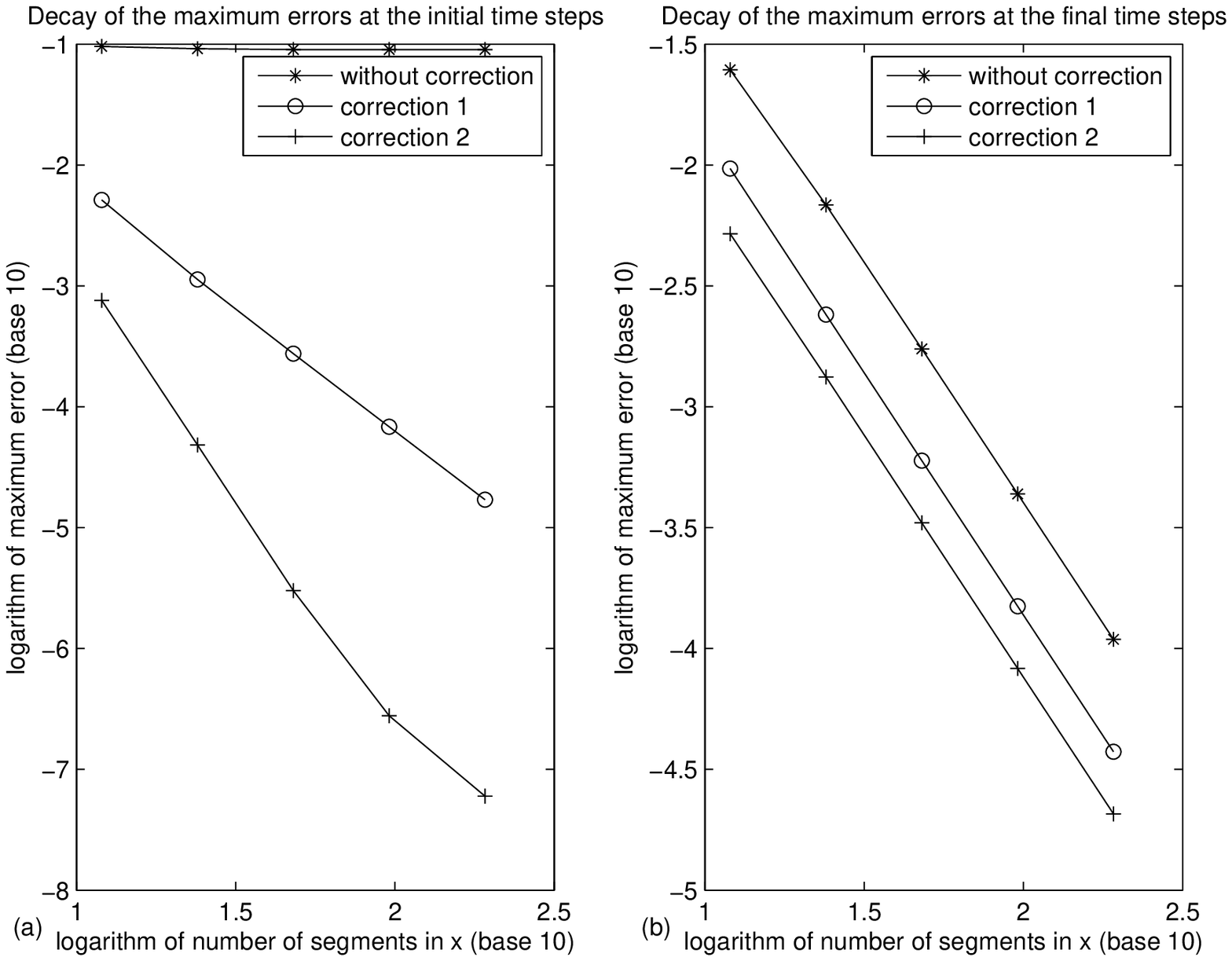}}
\caption{Decay of the maximum errors. (a) at initial step.(b)at final step.}
\label{f8}
\end{figure}
As in Section \ref{s3.2}, we take $\nu=0.2,\quad g_{1}=0,\quad g_{2}=0$ and $h=\sin(\frac{7\pi}{4}x+\frac{1}{4}\pi)$. It is easy to check that, for this test case, both the zeroth and the first order compatibility conditions at the right corner are met, but neither of them are met at the left corner.\\
\vspace{0.1cm}

Generally, given arbitrary initial and boundary conditions, it is impossible to find an analytic solution for the nonlinear reaction-diffusion equation. Therefore, in general, it is impossible to compute the real errors. As an alternative, we compute the comparative errors, which are the differences between the two numerical solutions for the problem, one with the mesh under consideration, and the other one with a finer mesh. In the following the term error is to be understood in this sense.\\
\vspace{0.1cm}
We first compute the solution of \eqref{e3.1} without any
correction procedure, i.e.~with $S=0$ in \eqref{e4.9}.
The solution is plotted in Fig.~\ref{f5} (a), and it
displays a sharp gradient around the left corner
of the time--space axes due to the incompatibility
between the initial and boundary conditions there.
In order to have an overview
of the structure of the errors in the solution we plot
the pointwise errors in Fig.~\ref{f5} (b). The pike appears
near the left corner of the time--space axes, as expected,
and it dissipates away quickly.
For comparison we plot, in Fig.~\ref{f6},
the solution and the pointwise
errors computed with {\it Correction Procedure 1},
and in Fig.~\ref{f9},
the solution and the pointwise
errors computed with {\it Correction Procedure 2}.
With {\it Correction procedure 1} the magnitude of
the errors at the left corner (see Fig.~\ref{f6} (b))
has been reduced by roughly two orders
(compared to Fig.~\ref{f5} (b))! The most pronounced
errors actually appear later in the simulation due to
accumulation of the errors. With
{\it Correction procedure 2} the errors
(Fig.~\ref{f9}) are
further reduced by roughly one half
(compared to Fig.~\ref{f6}).

The evolution of the maximum errors at each time step
is plotted in Fig.~\ref{f7} (a).
The evolution of the maximum errors, with
{\it Correction Procedure 1}
 enabled, is displayed in Fig.~\ref{f7} (b).
The magnitude of the maximum errors has been
reduced by indeed two orders
(compared to Fig.~\ref{f7} (a)).

When \emph{Correction Procedure 2} is enabled,
we see further improvement in the accuracy, though
it is less dramatic. For comparison we plot the result
in the same figure (Fig.~\ref{f7} (b)) as that for the result with
\emph{Correction Procedure 1}, and we see that the
magnitude of the maximum errors is roughly halved.

\subsection{Comparison of the convergence rates}\label{s4.3}
In this subsection we study the decay of the maximum errors as functions of the mesh size(number of the segments in $x$),
with and without the correction procedures applied.
The results confirm the effectiveness of
the correction procedures.
Fig.~\ref{f8} (b) shows that, with and without the
correction procedure applied,  the maximum errors
at a fixed time $t=0.05$ decay at approximately
the second order (the slope of the line),
which is the order of accuracy of the finite element
scheme. However, the maximum errors of the simulation
when {\it Correction procedure 1} is applied is
about one half order smaller in magnitude than without any correction procedure.
The maximum errors of the simulation with
{\it Correction procedure 2}
applied is smaller by an additional factor of about $0.25$.

The most interesting and informative comparison
can be made between the decay rates of the maximum
errors at the initial time step ($t=\Delta t$, $\Delta t$
varying with different configurations).
In Fig.~\ref{f8} (a),
the curve for the simulation without any correction
procedure stick to the upper frame of the figure;
the maximum errors will not come down whatsoever.
With {\it Correction procedure 1}, the maximum errors
decay at roughly the first order with respect
to $\Delta x$, and with {\it Correction procedure 2}, the
results are slightly better in terms of both the magnitude
of the maximum errors and the decay rate.

\section{Conclusion}\label{s4}
Incompatibilities between the initial and boundary conditions
for evolution PDEs have an adverse effect on the accuracy of numerical
simulations, especially near the time--space corners. No matter
how fine the grid is, the magnitude of the maximum
errors persists, with the pike of the errors moving towards the corners
as the grid gets finer.

For the same configuration, {\it Correction procedure 1} reduces the
magnitude of the maximum errors by about two orders;
and {\it Correction procedure 2} further reduces the magnitude
by roughly one half.

For a fixed time, the correction procedures have no effect
on the convergence rate of the solution at that point, but the
correction procedures always give more accurate results, with
{\it Correction procedure 2} being more effective than
{\it Correction procedure 1}. At the initial time step ($t=\Delta t$),
without any correction procedure the magnitude of the maximum
errors persists as mesh size gets small; with {\it Correction procedures
1 or 2}, the errors diminish at roughly order 1 with respect to
$\Delta x$, with Correction procedure 2 doing slightly better
than {\it Correction procedure 1}.

The approach described in this article for deriving correction procedures
does not
depend  on any particular property of the viscous Burgers equation
or the reaction--diffusion equation
other than its diffusiveness. Hence we believe that
these procedures can also be applied to other nonlinear
diffusive equations in space dimension one.
As mentioned earlier, in a work in progress \cite{ChQinTe2}
we study a totally different method related for higher space
dimensions.

\section*{Acknowledgments}
This work was supported in part by
NSF grants
DMS0604235 and DMS0906440 and by the Research Fund of Indiana University.

\bibliographystyle{amsplain}
\bibliography{biblio}

\end{document}